\documentclass{article}

\usepackage{amsmath}
\usepackage{url}

\begin{document}

\title{Evaluations For $\varsigma (2)$, $\varsigma (4)$, \ldots , $\varsigma (2k)$
Based On The WZ Method}
\author{ YiJun Chen}
\date{}
\maketitle

\begin{abstract}
Based on the framework of the WZ theory , a new evaluation for
$\varsigma (2) = \frac{\pi ^2}{6}$ and $\varsigma (4) = \frac{\pi ^4}{90}$
was given respectively, finally, a recurrence formula for
$\varsigma (2k)$, which is equivalent to the classical formula $B_{2k}(\frac{1}{2})=(2^{-2k+1}-1)B_{2k}$, was given.
\end{abstract}

\section{Introduction, Lemmas and Main Results.}

We know that there are many evaluations (or proofs) for $\varsigma (2) =
\frac{\pi ^2}{6}$ since the first evaluation belong to Euler, e.g. see \cite{chapman}-\cite{srivastava}
and the related references therein. We also know that there are two recurrence formulas for $\varsigma (2k)$ in \cite{srivastava}

\noindent
$
\varsigma (2n)= ( - 1)^{n - 1}\frac{(2\pi )^{2n - 1}}{(2n)!(2^{2n - 1} -
1)}\left[ {\frac{\pi }{2(2n + 1)} + \sum\limits_{j = 1}^{n - 1} {( -
1)^j{{2n} \choose {2j - 1}}\frac{(2j - 1)!(2^{2j - 1} - 1)}{(2\pi )^{2j - 1}}}
\varsigma (2j)} \right],
$

\noindent
$
\varsigma (2n) = ( - 1)^{n - 1}\frac{(2\pi )^{2n - 1}}{(2n - 1)!(2^{2n} -
1)}\left[ {\frac{\pi }{4n} + \sum\limits_{j = 1}^{n - 1} {( - 1)^j{{2n - 1} \choose {2j - 1}}\frac{(2j - 1)!(2^{2j - 1} - 1)}{(2\pi )^{2j - 1}}}
\varsigma (2j)} \right].
$

In this paper, I will give a new evaluation for $\varsigma (2) = \frac{\pi
^2}{6}$ based on the framework of WZ theory (see \cite{wilf}-\cite{zeilberger}), repeating the
process of evaluation can also be applied to evaluating $\varsigma (4)$, $\varsigma (6)$, \ldots.
Finally, through the same process repeatedly, I obtained a recurrence formula for
$\varsigma (2k)$ which is similar to, but different from the recurrence formulas for
$\varsigma (2k)$ mentioned above. As
$\varsigma (2k) = \frac{2^{2k - 1}( - 1)^{k - 1}B_{2k} }{(2k)!}\pi ^{2k}$
(belong to Euler too), and $B_k $(where $k \in N_0 = N \cup \{0\} $ ) is called $k$-th Bernoulli
number, thus by the recurrence formula for $\varsigma (2k)$ in the following theorem, we can obtain a formula for Bernoulli polynomial $B_{2k}$: $B_{2k}(\frac{1}{2})=(2^{-2k+1}-1)B_{2k}$, where $B_{n}(x)$ is the Bernoulli polynomial of order $n$, in fact, there are equivalent.
The following theorem is the main result in this paper.

\noindent
\textbf{Theorem. }Given $\varsigma (s) = \sum\limits_{n = 1}^{ + \infty }
{\frac{1}{n^s}} $ (where $Re(s) > 1)$, then we have $\varsigma (\ref{eq2}) =
\frac{\pi ^2}{6}$, $\varsigma (4) = \frac{\pi ^4}{90}$, more generally with
the convention $\sum\limits_{k = 1}^0 {a(k)} = 0$, for $\varsigma (2l)$
(where $l \in N )$ the following recurrence formula hold
\[
\varsigma (2l) = \left( {\frac{2^{2l - 1}}{1 - 2^{2l}}} \right)\left\{
{\left[ {\frac{( - 1)^{l + 1}}{4l} + \frac{( - 1)^l}{2}} \right]\frac{\pi
^{2l}}{\Gamma (2l)} + \sum\limits_{j = 1}^{l - 1} {\frac{( - 1)^{l - j}\pi
^{2(l - j)}}{\Gamma (2(l - j) + 1)}\varsigma (2j)} } \right\}
\]
\noindent
where $\Gamma (z)$ is gamma function.

To prove the theorem, we need the following lemmas.

\noindent
\textbf{Lemma 1. }If given a continuous-discrete WZ pair $(F(x,k),G(x,k))$,
that is $(F(x,k),G(x,k))$ satisfy the following so called
continuous-discrete WZ equation

\begin{equation}
\label{eq1}
\frac{\partial F(x,k)}{\partial x} = G(x,k + 1) - G(x,k)
\end{equation}

\noindent
then for all $m, n \in N_0 $, for all $h, x \in R$, we
have

\[
\sum\limits_{k = m}^n {F(x,k)} - \sum\limits_{k = m}^n {F(h,k)} = \int_h^x
{G(t,n + 1)dt} - \int_h^x {G(t,m)dt} .
\]

\noindent
\textbf{Lemma 2. }If for $a, x \in R$, $a<x$, $f(x)$ is
integrable on the interval $(a,x)$, then we have

\[
\int_a^x {\left( {\int_a^{t_1 } { \cdots \left( {\int_a^{t_{k - 1} } {f(t_k
)} dt_k } \right)} \cdots dt_2 } \right)} dt_1 = \frac{1}{\Gamma
(k)}\int_a^x {(x - t)^{k - 1}f(t)} dt.
\]

\noindent
\textbf{Lemma 3. } For all $n \in N$, we have

\[
\sum\limits_{k = 1}^n {\cos kx} = - \frac{1}{2} + \frac{1}{2}\frac{\sin
\left[ {(2n + 1)x / 2} \right]}{\sin (x / 2)}.
\]

\noindent
\textbf{Lemma 4. } For all $n \in N_0 $, we have

\[
\int_0^\pi {\frac{\sin \left[ {(2n + 1)x / 2} \right]}{\sin (x / 2)}dx = \pi
} .
\]

\noindent
\textbf{Lemma 5. } For all $s \ge 1$, we have

\[
\mathop {\lim }\limits_{n \to + \infty } \int_0^\pi {\frac{x^s\sin \left[
{(2n + 1)x / 2} \right]}{\sin (x / 2)}dx = 0} .
\]

As the proof of \textbf{Lemma 1} is easy, we omit the details of proof here.
\textbf{Lemma 2} can be seen in \cite{melzak}-\cite{iwasaki}, \textbf{Lemma 3 }and \textbf{Lemma 4}
can be seen in \cite{byerly}. The proof of \textbf{Lemma 5} will be given below.

\section{Proof of Lemma 5.}

Let $t = \frac{x}{2}$, then we have

\[
\int_0^\pi {\frac{x^s\sin \left[ {(2n + 1)x \mathord{\left/ {\vphantom {x
2}} \right. \kern-\nulldelimiterspace} 2} \right]}{\sin \left( {x
\mathord{\left/ {\vphantom {x 2}} \right. \kern-\nulldelimiterspace} 2}
\right)}} dx = 2^{s + 1}\int_0^{\frac{\pi }{2}} {\frac{t^s\sin \left[ {(2n +
1)t} \right]}{\sin (t)}dt} .
\]
Let $f(t) = \left\{ {{\begin{array}{*{20}c}
 {t^s\csc t} \hfill & {0 < t < \textstyle{\pi \over 2}} \hfill \\
 0 \hfill & {t = 0} \hfill \\
\end{array} }} \right.$, it is an easy exercise of calculus that when $s \ge 1$, $f(t)$ is
differentiable and monotone (increasing) on $[0,\textstyle{\pi \over 2}]$,
then by \textbf{The Second Mean Value Theorem For Integrals}, we know that there exist $\xi$ on
$[0,\textstyle{\pi \over 2}]$ such that

\begin{eqnarray*}
\lefteqn{\left| {\int_0^{\frac{\pi }{2}} {\frac{t^s}{\sin t}\sin ((2n + 1)t)dt} }
\right|}\\
 &=& \left| {\int_0^{\frac{\pi }{2}} {f(t)\sin ((2n + 1)t)dt} }
\right|\\
 &=& \left| {f(0 + 0)\int_0^\xi {\sin ((2n + 1)t)dt}+f\left( {\frac{\pi }{2}- 0} \right)\int_\xi ^{\frac{\pi }{2}} {\sin ((2n + 1)t)dt} } \right|\\
 & = & \left( {\frac{\pi }{2}} \right)^s\left| {\frac{1}{2n + 1}\left. {( - \cos
((2n + 1)t))} \right|{\begin{array}{*{20}c}
 {\textstyle{\pi \over 2}} \hfill \\
 \xi \hfill \\
\end{array} }} \right|\\
 &\le& \left( {\frac{\pi }{2}} \right)^s\frac{2}{2n + 1}.
\end{eqnarray*}
By the result above, we can conclude that

\[
\mathop {\lim }\limits_{n \to + \infty } \int_0^\pi {\frac{x^s\sin \left(
{\frac{(2n + 1)x}{2}} \right)}{\sin \left( {\frac{x}{2}} \right)}} dx =
\mathop {\lim }\limits_{n \to + \infty } 2^{s + 1}\int_0^{\frac{\pi }{2}}
{\frac{t^s\sin ((2n + 1)t)}{\sin t}} dt = 0,
\]
the proof of \textbf{Lemma 5} was completed.

\noindent
\textbf{Remarks:} 1. It is worth mentioning that we can prove\textbf{ Lemma
5} by Riemann-Lebesgue lemma directly as follows. Let $f(t) = \left\{ {{\begin{array}{*{20}c}
 {\left( {\frac{t}{2}} \right)^s\csc \left( {\frac{t}{2}} \right)} \hfill &
{0 < t < \pi } \hfill \\
 0 \hfill & {t = 0} \hfill \\
\end{array} }} \right.$, it is an easy exercise of calculus that $f(t)$ is continuous on $[0,\pi]$, of course, $f(t)$ is Riemann integrable on $[0,\pi]$, then by Riemann-Lebesgue lemma, we have

\[
\mathop {\lim }\limits_{n \to + \infty } \int_0^\pi {f(x)\sin \left(
{\frac{(2n + 1)x}{2}} \right)} dx =0,
\]
because

\[
\int_0^\pi {f(x)\sin \left(
{\frac{(2n + 1)x}{2}} \right)} dx =\int_0^\pi {\frac{\left( {\frac{x}{2}} \right)^s\sin \left(
{\frac{(2n + 1)x}{2}} \right)}{\sin \left( {\frac{x}{2}} \right)}} dx,
\]
finally, we have

\[
\mathop {\lim }\limits_{n \to + \infty } \int_0^\pi {\frac{x^s\sin \left(
{\frac{(2n + 1)x}{2}} \right)}{\sin \left( {\frac{x}{2}} \right)}} dx =0.
\]

\noindent
2. It is also worth mentioning that
when $s \ge 2$, we can prove\textbf{ Lemma 5} by using integration by parts,
but when $1 \le s < 2$, the method can't be used.

\section{Proof of The Theorem.}

\textbf{(A) }\textbf{Proof of} $\varsigma (2) = \frac{\pi ^2}{6}.$
Setting $F_1 (x,k) = \frac{\cos (kx)}{k^2}$, $G_1 (x,k) = \sum\limits_{j =
1}^{k - 1} {\frac{ - \sin (jx)}{j}} $, then it is easy to
verify that $(F_1(x,k),G_1(x,k))$ is a continuous-discrete WZ pair, that is,
they satisfy the equation (\ref{eq1}). Now let $h = 0$, $m = 1$, with the
convention $\sum\limits_{k = 1}^0 {a(k)} = 0$, by using\textbf{ Lemma 1 }we
get

\begin{equation}
\label{eq2}
\sum\limits_{k = 1}^n {\frac{\cos (kx)}{k^2}} - \sum\limits_{k = 1}^n
{\frac{1}{k^2}} = \int_0^x {G_1 (t,n + 1)dt}
\end{equation}
To evaluate $G_1 (x,n + 1) = \sum\limits_{j = 1}^n {\frac{ - \sin
(jx)}{j}}$, we also use\textbf{ Lemma 1}. Now set

\[
F_2 (x,k) = \frac{ - \sin (kx)}{k},
\quad
G_2 (x,k) = \sum\limits_{j = 1}^{k - 1} { - \cos (kx)}  ,
\]
then it is easy to verify that $(F_2 (x,k),G_2 (x,k))$ is a continuous-discrete WZ
pair, and for all $k \in N$, the following result hold $F_2 (0,k) = 0$. With
the convention $\sum\limits_{k = 1}^0 {a(k)} = 0$, by using\textbf{
Lemma 1} and \textbf{Lemma 3}, we obtain

\begin{equation}
\label{eq3}
\sum\limits_{k = 1}^n {\frac{ - \sin (kx)}{k}} = \int_0^x {G_2 (t,k)} dt =
\int_0^x {\left\{ {\frac{1}{2} - \frac{\sin \left[ {(2n + 1)t / 2}
\right]}{2\sin (t / 2)}} \right\}} dt
\end{equation}
By using (\ref{eq2}), (\ref{eq3}) and\textbf{ Lemma 2}, we obtain

\begin{eqnarray*}
\sum\limits_{k = 1}^n {\frac{\cos (kx)}{k^2}} - \sum\limits_{k = 1}^n
{\frac{1}{k^2}} &=& \int_0^x {G_1 (t,n + 1)} dt\\
 &=& \int_0^x {\left\{ {\int_0^{t_1 } {\left[ {\frac{1}{2} - \frac{\sin \left[
{(2n + 1)t_2 / 2} \right]}{2\sin (t_2 / 2)}} \right]dt_2 } } \right\}dt_1 }\\
&=&\frac{1}{2}\int_0^x {(x - t)} dt - \frac{x}{2}\int_0^x \frac{\sin \left[
{(2n + 1)t / 2} \right]}{\sin (t / 2)}dt\\
&&+\:\frac{1}{2}\int_0^x {\frac{t\sin\left[ {(2n + 1)t / 2} \right]}{\sin (t / 2)}dt}\\
&=&I_1 (x) + I_2 (x) + I_3 (x).
\end{eqnarray*}
Recalling $\sum\limits_{k = 1}^{ + \infty } {\frac{( - 1)^k}{k^s}} -
\varsigma (s) = \left( { - 2 + \frac{1}{2^{s - 1}}} \right)\varsigma (s)$,
let $x = \pi $ at first, and then let $n \to + \infty $, we conclude that

\[
\mathop {\lim }\limits_{n \to + \infty } \left[ {\sum\limits_{k = 1}^n
{\frac{\cos (k\pi )}{k^2}} - \sum\limits_{k = 1}^n {\frac{1}{k^2}} } \right]
= \sum\limits_{k = 1}^{ + \infty } {\frac{( - 1)^k}{k^2}} - \varsigma (\ref{eq2}) =
- \frac{3}{2}\varsigma (\ref{eq2}).
\]
After some computations, we obtain

\[
I_1 (\pi ) = \frac{1}{2}\int_0^\pi {(\pi - t)} dt = \frac{\pi ^2}{4}.
\]
By \textbf{Lemma 4}, we obtain

\[
I_2 (\pi ) = - \frac{\pi }{2}\int_0^\pi {\frac{\sin \left[ {(2n + 1)t / 2}
\right]}{\sin (t / 2)}} dt = - \frac{\pi ^2}{2}.
\]
Now by\textbf{ Lemma 5}, we obtain

\[
\mathop {\lim }\limits_{n \to + \infty } I_3 (\pi ) = \mathop {\lim
}\limits_{n \to + \infty } \frac{1}{2}\int_0^\pi {\frac{t\sin \left[ {(2n +
1)t / 2} \right]}{\sin (t / 2)}dt = 0} .
\]
Finally we conclude that $ - \frac{3}{2}\varsigma (\ref{eq2}) = \frac{\pi ^2}{4} -
\frac{\pi ^2}{2} = - \frac{\pi ^2}{4}$, that is $\varsigma (2) = \frac{\pi
^2}{6}$\textbf{.}

\noindent
\textbf{(B) }\textbf{Proof of }$\varsigma (4) = \frac{\pi ^4}{90}.$ Setting

\[
F_1 (x,k) = \frac{\cos (kx)}{k^4},
\quad
G_1 (x,k) = \sum\limits_{j = 1}^{k - 1} { \frac{-\sin (jx)}{j^3}} ,
\]

\[
F_2 (x,k) = \frac{ - \sin (kx)}{k^3},
\quad
G_2 (x,k) = \sum\limits_{j = 1}^{k - 1} { \frac{-\cos (jx)}{j^2}},
\]

\[
F_3 (x,k) =  \frac{-\cos (kx)}{k^2},
\quad
G_3 (x,k) = \sum\limits_{j = 1}^{k - 1} {\frac{\sin (jx)}{j}} ,
\]

\[
F_4 (x,k) = \frac{\sin (kx)}{k},
\quad
G_4 (x,k) = \sum\limits_{j = 1}^{k - 1} {\cos (jx)} ,
\]
It is easy to verify that for $j = 1,2,3,4$, $(F_j (x,k),G_j (x,k))$ satisfy
equation (\ref{eq1}). Setting $H_n^{(l)} = \sum\limits_{k = 1}^n {\frac{1}{k^l}} $,
completely analogous to the proof of $\varsigma (\ref{eq2}) = \frac{\pi ^2}{6}$( some details are omitted here. ), we get

\begin{eqnarray*}
\sum\limits_{k = 1}^n {\frac{\cos (kx)}{k^4}} - H_n^{(4)} &=&\frac{1}{\Gamma
(4)}\left\{ {\left[ {\int_0^x { - \frac{1}{2}(x - t)^3dt +
\frac{x^3}{2}\int_0^x {\frac{\sin \left[ {(2n + 1)t / 2} \right]}{\sin (t /
2)}} dt} } \right]} \right\}\\
&&+\:\frac{1}{\Gamma (4)}\left\{ {\sum\limits_{k = 1}^2 {\left(
{{\begin{array}{*{20}c}
 3 \hfill \\
 k \hfill \\
\end{array} }} \right)\frac{x^k}{2}\int_0^x {( - t)^{3 - k}\frac{\sin \left[
{(2n + 1)t / 2} \right]}{\sin (t / 2)}} dt} } \right\}\\
 &&-\: \frac{1}{\Gamma (4)}\int_0^x {(x - t)H_n^{(\ref{eq2})} } dt\\
&=&I_1 (x) + I_2 (x) + I_3 (x).
\end{eqnarray*}

Setting $x = \pi $, by \textbf{Lemma 4} we obtain $I_1 (\pi ) = \frac{\pi
^4}{16}$, and by\textbf{ Lemma 5} we obtain $\mathop {\lim }\limits_{n \to +
\infty } I_2 (\pi ) = 0$, recall$\mathop {\lim }\limits_{n \to + \infty }
H_n^{(\ref{eq2})} = \varsigma (\ref{eq2}) = \frac{\pi ^2}{6}$, we obtain $\mathop {\lim
}\limits_{n \to + \infty } I_3 (\pi ) = - \frac{\pi ^4}{12}$. Recall
$\mathop {\lim }\limits_{n \to + \infty } H_n^{(4)} = \varsigma (4)$, we get

\[
\mathop {\lim }\limits_{n \to + \infty } \left[ {\sum\limits_{k = 1}^n
{\frac{\cos (k\pi )}{k^4}} - H_n^{(4)} } \right] = \sum\limits_{k = 1}^{ +
\infty } {\frac{( - 1)^k}{k^4}} - \varsigma (4) = \left( { - 2 +
\frac{1}{2^3}} \right)\varsigma (4).
\]

Finally we obtain $\varsigma (4) = \frac{\pi ^4}{90}$.

\noindent
\textbf{(C) }\textbf{Proof of }$\varsigma (2l) = \left( {\frac{2^{2l - 1}}{1
- 2^{2l}}} \right)\left\{ {\left[ {\frac{( - 1)^{l + 1}}{4l} + \frac{( -
1)^l}{2}} \right]\frac{\pi ^{2l}}{\Gamma (2l)} + \sum\limits_{j = 1}^{l - 1}
{\frac{( - 1)^{l - j}\pi ^{2(l - j)}}{\Gamma (2(l - j) + 1)}\varsigma (2j)}
} \right\}.$ The result can be proved in the above framework of proving $\varsigma (\ref{eq2}) =
\frac{\pi ^2}{6}$ and $\varsigma (4) = \frac{\pi ^4}{90}$, some details are omitted here. For convenience, setting
\[
H_n^{(l)} (x) = \sum\limits_{k = 1}^n {\frac{\cos (kx)}{k^l}} ,
\quad
H_n^{(l)} (0) = H_n^{(l)} ,
\quad
H_n^{(l)} (\pi ) = \sum\limits_{k = 1}^n {\frac{( - 1)^k}{k^l}} ,
\]

\[
I_j (f)(x) = \frac{1}{\Gamma (j)}\int_0^x {(x - t)^{j - 1}f(t)dt} ,
\]
with the convention that $I_0 (f)(x) = f(x)$, where $j \in N_0 $. Now it is
easy to verify that for all $j \in N_0 $, $I_j $ own the following
properties

\[
I_j (f + g)(x) = I_j (f)(x) + I_j (g)(x),
\quad
I_j (cf)(x) = cI_j (f)(x),
\]
where $c$ is a constant having nothing to do with $t$, the variable of
integral. Also setting

\[
f(t) = - \frac{1}{2} + \frac{1}{2}\frac{\sin \left[ {(2n + 1)t / 2}
\right]}{\sin (t / 2)},
\]
we obtain

\[
H_n^{(2l)} (x) = ( - 1)^lI_{2l} (f)(x) + \sum\limits_{j = 1}^l {( - 1)^{l -
j}I_{2(l - j)} (H_n^{(2j)} )(x)} .
\]
Let $x = \pi $, then we get the following result

\begin{eqnarray*}
\sum\limits_{k = 1}^n {\frac{( - 1)^k}{k^{2l}}} - H_n^{(2l)} &=& ( -
1)^lI_{2l} \left( { - \frac{1}{2}} \right)(\pi ) + ( - 1)^lI_{2l} \left(
{\frac{\sin \left[ {(2n + 1)t / 2} \right]}{2\sin (t / 2)}} \right)(\pi )\\
&&+\: \sum\limits_{j = 1}^{l - 1} {( - 1)^{l - j}I_{2(l - j)} (H_n^{(2j)} )(\pi
)}\\
&=& \textrm{I} + \textrm{II} + \textrm{III}.
\end{eqnarray*}
Next, let us consider $\textrm{I}$, $\textrm{II}$ and $\textrm{III}$ respectively

\begin{eqnarray*}
\textrm{I} &=& ( - 1)^l\left( { - \frac{1}{2}} \right)\frac{1}{\Gamma (2l)}\int_0^\pi
{(\pi - t)^{2l - 1}dt = \frac{( - 1)^{l + 1}}{4l}} \frac{\pi ^{2l}}{\Gamma
(2l)}\\
\textrm{II} &=& ( - 1)^l\frac{1}{2\Gamma (2l)}\int_0^\pi {(\pi - t)^{2l - 1}\frac{\sin
\left[ {(2n + 1)t / 2} \right]}{\sin (t / 2)}} dt\\
 &=& ( - 1)^l\frac{1}{2\Gamma (2l)}\int_0^\pi {\pi ^{2l - 1}\frac{\sin \left[
{(2n + 1)t / 2} \right]}{\sin (t / 2)}dt}\\
 &&+\: ( - 1)^l\frac{1}{2\Gamma (2l)}\int_0^\pi {\sum\limits_{k = 1}^{2l - 2}
\binom{2l-1}{k}\pi ^k( - t)^{2l - 1 - k} \frac{\sin \left[ {(2n +
1)t / 2} \right]}{\sin (t / 2)}dt}\\
&=& \textrm{II}_1 + \textrm{II}_2 .
\end{eqnarray*}
By Lemma 4, we obtain

\[
\textrm{II}_1 = ( - 1)^l\frac{\pi ^{2l}}{2\Gamma (l)},
\]
and by Lemma 5 we conclude that $\mathop {\lim }\limits_{n \to + \infty }
\textrm{II}_2 = 0$. By using the results above, we obtain

\[
\mathop {\lim }\limits_{n \to + \infty } \textrm{II} = ( - 1)^l\frac{\pi
^{2l}}{2\Gamma (l)}.
\]
After some computations, we obtain

\begin{eqnarray*}
\textrm{III} &= &\sum\limits_{j = 1}^{l - 1} {( - 1)^{l - j}H_n^{(2j)} \frac{1}{\Gamma
(2(l - j))}\int_0^\pi {(\pi - t)^{2(l - j) - 1}dt} }\\
 &= &\sum\limits_{j =
1}^{l - 1} {( - 1)^{l - j}H_n^{(2j)} \frac{\pi ^{2(l - j)}}{\Gamma \left(
{2(l - j) + 1} \right)}} .
\end{eqnarray*}
Recalling$\mathop {\lim }\limits_{n \to + \infty } H_n^{(2j)} = \varsigma
(2j)$, we conclude that

\[
\mathop {\lim }\limits_{n \to + \infty } \textrm{III} =
\sum\limits_{j = 1}^{l - 1} {\frac{( - 1)^{l - j}\pi ^{2(l - j)}}{\Gamma
(2(l - j) + 1)}\varsigma (2j)}.
\]
It is easy to verify that

\[
\mathop {\lim }\limits_{n \to + \infty } \left( {\sum\limits_{k = 1}^n
{\frac{( - 1)^k}{k^{2l}}} - \sum\limits_{k = 1}^n {\frac{1}{k^{2l}}} }
\right) = \left( { - 2 + \frac{1}{2^{2l - 1}}} \right)\varsigma (2l).
\]
So finally with the convention $\sum\limits_{k = 1}^0 {a(k)} = 0$, we obtain
the following recurrence formula for $\varsigma (2k)$

\[
\varsigma (2l) = \left( {\frac{2^{2l - 1}}{1 - 2^{2l}}} \right)\left\{
{\left[ {\frac{( - 1)^{l + 1}}{4l} + \frac{( - 1)^l}{2}} \right]\frac{\pi
^{2l}}{\Gamma (2l)} + \sum\limits_{j = 1}^{l - 1} {\frac{( - 1)^{l - j}\pi
^{2(l - j)}}{\Gamma (2(l - j) + 1)}\varsigma (2j)} } \right\}.
\]
The proof of the theorem is completed.

\noindent
\textbf{Remarks: }1. We can set

\[
F_1 (x,k) = \frac{e^{ikx}}{k^2},
\quad
G_1 (x,k) = \sum\limits_{j = 1}^{k - 1} {\frac{ie^{ijx}}{j}} ,
\]

\[
F_2 (x,k) = \frac{ie^{ikx}}{k},
\quad
G_2 (x,k) = \sum\limits_{j = 1}^{k - 1} { - e^{ijx}}
\]
where $i = \sqrt { - 1} $, it is easy to verify that $(F_j (x,k),G_j (x,k))$
(where $j = 1,2)$ is a continuous-discrete WZ pair. Then through the same
process above, we can also obtained $\varsigma (\ref{eq2}) = \frac{\pi ^2}{6}$, the
details are omitted here. Of course, we can do in the same way for $\varsigma (4)$, $\varsigma (6)$, \ldots, and for the general case $\varsigma (2k)$ respectively. 2. It is also worth mentioning that the ideas in the proof above can be used to solve other similar problems of summation of infinite series, and I will give the details in another paper.

\bigskip

\noindent\textit{School of Mathematical Science,
South China Normal University, Guangzhou, P.R.China\\
chenyijun73@yahoo.com.cn}

\bigskip


\begin{thebibliography}{9}
\bibitem{chapman} Robin Chapman, Evaluating $\varsigma (\ref{eq2})$, available at \url{http://empslocal.ex.ac.uk/people/staff/rjchapma/etc/zeta2.pdf}
\bibitem{kalman1} Dan Kalman, Six ways to sum a series, College Math. J. 24(1993), 402-421.
\bibitem{kalman2} Dan Kalman, Mark McKinzie, Another way to sum a series-historical appendix,
    details of an historical investigation as to whether Euler knew Lewin's argument(2011), available at \url{http://www.dankalman.net/eulerdilog/.}
\bibitem{kalman3} Dan Kalman, Mark McKinzie, Another way to sum a series: Generating Functions, Euler, and The Dilog Function, Amer. Math.Monthly 119(2012), 42-51.
\bibitem{srivastava} H. M. Srivastava, Junesang Choi, Zeta and q-Zeta Functions and Associated Series and Integrals, New York: Elsevier Inc. , 2012
\bibitem{wilf} H. S. Wilf, D. Zeilberger, Rational functions certify combinatorial identities, J. of Amer. Math. Soc. 3(1990), 147-158.
\bibitem{petkov} M. Petkov\v{s}ek, H. S. Wilf, D.Zeilberger, A=B, MA: A.K.Peters Wellesly, 1996
\bibitem{zeilberger} D. Zeilberger, Closed form (Pun intended!), Contemp. Math. 143(1993), 579-607.
\bibitem{melzak} Z. A. Melzak, Companion to Concrete Mathematics: Mathematical Techniques and Various Applications, New York: John Wiley{\&}Sons, Inc. ,1973
\bibitem{iwasaki}K. Iwasaki, H. Kimura, Sh. Shimomura, M. Yoshida, From Gauss to Painlev\'{e}: A Modern Theory of Special Functions, Braunschweig: Vieweg, 1991
\bibitem{byerly}William Elwood Byerly, An Elementary Treatise on Fourier's Series and Spherical, Cylindrical, and Ellipsoidal Harmonics, With Applications to Problems in Mathematical Physics, New York: Ginn{\&}Company, 1893
\end{thebibliography}
\end{document}